\documentstyle[12pt]{article}

\begin{document}


\newtheorem{guess}{Theorem}[section]
\newcommand{\bth}{\begin{guess}$\!\!\!${\bf .}~}
\newcommand{\eeth}{\end{guess}}

\newtheorem{propo}[guess]{Proposition}
\newcommand{\bprop}{\begin{propo}$\!\!\!${\bf .}~}
\newcommand{\eprop}{\end{propo}}

\newtheorem{lema}[guess]{Lemma}
\newcommand{\blem}{\begin{lema}$\!\!\!${\bf .}~}
\newcommand{\elem}{\end{lema}}

\newtheorem{defe}[guess]{Definition}
\newcommand{\bdefe}{\begin{defe}$\!\!\!${\bf .}~}
\newcommand{\edefe}{\end{defe}}

\newtheorem{coro}[guess]{Corollary}
\newcommand{\bcor}{\begin{coro}$\!\!\!${\bf .}~}
\newcommand{\ecor}{\end{coro}}

\newtheorem{rema}[guess]{Remark}
\newcommand{\brem}{\begin{rema}$\!\!\!${\bf .}~\rm}
\newcommand{\erem}{\end{rema}}

\newcommand{\sing}{{\rm Sing}\,}
\newcommand{\lr}{\longrightarrow}
\newcommand{\ce}{{\cal E}}
\newcommand{\co}{{\cal O}}
\newcommand{\Co}{{\cal O}}
\newcommand{\ct}{{\cal T}}
\newcommand{\cf}{{\cal F}}
\newcommand{\cg}{{\cal G}}
\newcommand{\cv}{{\cal V}}
\newcommand{\ci}{{\cal I}}
\newcommand{\cn}{{\cal N}}
\newcommand{\cu}{{\cal U}}
\newcommand{\ch}{{\cal H}}
\newcommand{\cm}{{\cal M}}
\newcommand{\cc}{{\cal C}}
\newcommand{\cw}{{\cal W}}
\newcommand{\cp}{{\cal P}}
\newcommand{\cl}{{\cal L}}
\newcommand{\fq}{{q}}
\newcommand{\cd}{{\cal D}}
\newcommand{\Sing}{{\rm Sing}\,}
\newcommand{\eee}{{\rm End}\,}
\newcommand{\rank}{{\rm rank}\,}
\newcommand{\codim}{{\rm codim}\,}
\newcommand{\coker}{{\rm coker}\,}
\newcommand{\Aut}{{\rm Aut}\,}
\newcommand{\Ext}{{\rm Ext}\,}
\newcommand{\ad}{{\rm ad}\,}
\newcommand{\lm}{{\mbox{\Large $|$}}}
\newcommand{\rk}{{\rm rank}\,}

\newcommand{\ctext}[1]{\makebox(0,0){#1}}
\newcommand{\Bbb}{\bf}
\setlength{\unitlength}{0.1mm}

\title{Stability of the Poincar\'e Bundle}
\author{V. Balaji, L. Brambila Paz and P. E.
Newstead\thanks{All the authors are members of the VBAC
Research group of Europroj.  The work forms part of a
project ``Moduli of bundles in algebraic geometry'', funded
by the EU International Scientific Cooperation Initiative
(Contract no. CI1$^*$-CT93-0031).  The first author also
acknowledges support from the EPSRC (Grant no. GR/H38878)
and the second author from CONACYT (Project no.
3231-E9307).}}\date{30 May 1995}
\maketitle

\begin{abstract}
Let $C$ be a nonsingular projective curve of genus
$g\ge2$ defined over the complex numbers, and let
$M_{\xi}$ denote the moduli space of stable bundles of
rank $n$ and determinant $\xi$ on $C$, where $\xi$ is a
line bundle of degree $d$ on $C$ and $n$ and $d$ are
coprime. It is shown that the universal bundle
$\cu_{\xi}$ on $C\times M_{\xi}$ is stable with respect
to any polarisation on $C\times M_{\xi}$. It is shown
further that the connected component of the moduli
space of $\cu_{\xi}$ containing $\cu_{\xi}$ is
isomorphic to the Jacobian of $C$.
\end{abstract}
\pagebreak
\section*{Introduction}

In the study of moduli spaces of stable bundles on an
algebraic curve $C$, various bundles on the moduli space or
on the product of the moduli space with $C$ arise in a
natural way.  An interesting question to ask about any such
bundle is whether it is itself stable in some sense.

More precisely, let $C$ be a nonsingular projective curve of
genus $g\ge2$ defined over the complex numbers, and let $M
= M_{n,d}$ denote the moduli space of stable bundles of
rank $n$ and degree $d$ on
$C$, where $n$ and $d$ are coprime.  For any line bundle
$\xi$ of $\deg d$ on $C$, let $M_{\xi}$ denote the
subvariety of $M$ corresponding to bundles with
determinant $\xi$.  There exists on $C \times M$ a
universal (or Poincar\'e) bundle $\cu$ such that $\cu
\lm_{C \times \{ m\}}$ is the bundle on $C$ corresponding
to $m$.  Moreover the bundle $\cu$ is determined up to
tensoring with a line bundle lifted from $M$.

The direct image of $\cu$ on $M$ is called the Picard
sheaf of $\cu$; for $d> n (2g-2)$, this sheaf is a
bundle.   It was shown recently by Y.~Li [Li] that, if
$d>2gn$, this bundle is stable with respect to the ample
line bundle corresponding to the generalised theta divisor
(cf. [DN]). (Recall here that, if $H$ is an ample divisor
on a projective variety $X$, the {\it degree} $\deg E$ of a
torsion-free sheaf $E$ on $X$ is defined to be the
intersection number $[c_1(E) \cdot H^{\dim X-1}]$.  $E$ is
said to be {\it stable} with respect to $H$ (or
$H$-{\it stable}) if, for every proper subsheaf $F$ of $E$,
$$ \frac{\deg F}{\rk F} < \frac{\deg E}{\rk E}.$$  The
definition depends only on the polarisation defined by
$H$.) This extends previously known results for the case
$n=1$ ([U, Ke1, EL]).  We remark that the question of
stability of the Picard sheaf of $\cu_{\xi}$ is still open
(cf. [BV]).

In this paper, we investigate the stability of the
Poincar\'e bundle $\cu$ and its restriction $\cu_{\xi}$
to $M_{\xi}$ using methods similar to those of [Li].

Our main results are

\vspace{8pt}
\noindent
{\bf Theorem 1.5.}
$\cu_{\xi}$ is stable with respect to any polarisation on $C
\times M_{\xi}$.

\vspace{8pt}
\noindent
{\bf Theorem 1.6.}
$\cu$ is stable with respect to any polarisation of the form
$$a \alpha + b \Theta ,~~ a,b > 0$$
where $\alpha$ is ample on $C$ and $\Theta$ is the
generalised theta divisor on $M$.  (Note that $C$ has a
unique polarisation whereas $M$ does not.)

These results are proved in \S 1.  In \S 2 we discuss some
properties of the bundles $\eee \cu_{\xi}$ and $\ad
\cu_{\xi}$.  It is reasonable to conjecture that $\ad
\cu_{\xi}$ is also stable but we are able to prove this
only in the case $n=2$.

Finally, in \S 3 we consider the deformation theory of
$\cu_{\xi}$ using the results in \S 2.  The main result we
prove is that the only deformations of $\cu_{\xi}$ are
those of the form $\cu_{\xi} \otimes p_C^* L$, where $L$
is a line bundle of degree 0 on $C$.  More precisely, let
$H$ be any ample divisor on $C \times M_{\xi}$ and let
$M(\cu_{\xi})$ denote the moduli space of $H$-stable
bundles with the same numerical invariants as $\cu_{\xi}$
on $C \times M_{\xi}$; then

\vspace{8pt}
\noindent
{\bf Theorem 3.1.}
The connected component $M(\cu_{\xi})_0$ of $M(\cu_{\xi})$
containing \{$\cu_{\xi}$\} is isomorphic to the Jacobian
$J(C)$, the isomorphism $J(C) \lr M(\cu_{\xi})_0$ being
given by $$L \longmapsto \cu_{\xi} \otimes p_C^* L,$$
where $p_C: C \times M_{\xi} \lr C$ is the projection.

\vspace{8pt}
\noindent
{\bf Acknowledgement.}
Most
of the work for this paper was carried out during a visit by
the first two authors to Liverpool. They wish to
acknowledge the generous hospitality of the University of
Liverpool.

\renewcommand{\thesection}{\S \arabic{section}}
\section{Stability of $\cal U$}
\renewcommand{\thesection}{\arabic{section}}

We begin with some lemmas which are probably well known, but
which we could not find in the literature.

\blem \label{l11}
Let $X$ and $Y$ be smooth projective varieties of the same dimension
$m$.  Let $f: X - - \rightarrow Y$ be a dominant rational map defined
outside a subset $Z \subset X$ with $\codim_XZ\geq 2$.  Suppose that
$D_X$ and $D_Y$ are ample divisors on $X$ and $Y$ such that
$f^*D_Y|_{X-Z} \simeq D_X |_{X-Z}$.

Let $E$ be a vector bundle on $Y$ such that $f^* E$ extends to a vector
bundle $F$ on $X$. If $F$ is $D_X$-semi-stable (resp. stable) on
$X$, then $E$ is $D_Y$-semi-stable (resp. stable) on $Y$.

\elem

\noindent
{\bf Proof.}  The proof is fairly straight forward (cf. [Li,
pp.548, 549]).
Let rank $E = n$.  Suppose that $V$ is a torsion-free quotient
of $E$,
$$
E \lr V \lr 0 \eqno (1)
$$
with rank $V = r < n$. When $F$ is $D_X$-semi-stable, we
need to check the inequality:
$$
\frac{\deg V}{r} \geq \frac{\deg E}{n}.
\eqno (2)
$$
{}From (1), we have
$$
f^* E \lr f^* V \lr 0.
$$
If $G = f^*V/(torsion)$, then
$$
\deg f^* V \geq \deg G. \eqno (3)
$$
(Note that $\deg f^* V = [c_1(f^*V)\cdot D_X^{m-1}]$, where
$c_1(f^*V)$ makes sense since $f^*V$ is defined outside a
subset of codimension 2.)

Moreover we have
$$ 0 \lr G^* \lr f^* (E)^* $$
on X-Z.   Let $\widetilde{G}$ be a subsheaf of $F^*$
extending $G^*$.
Since $\codim_X Z \geq 2$, we have $\deg \widetilde{G} =
\deg G^* = - \deg G$.

By semi-stability of $F^*$,
$$
\frac{\deg\widetilde{G}}{r} \leq \frac{\deg F^*}{n}. \eqno
(4) $$
Note that
$$
\begin{array}{rcl}
\deg E & = & [c_1(E) . D_Y^{m-1}] \\ [2mm]
& = & [c_1 (f^*E) . D_X^{m-1}] (\deg f)^{-1} \\ [2mm]
& = & (\deg f^* E) (\deg f)^{-1}.\\
\end{array}
$$

The same applies to $V$, so we have
$$
\begin{array}{rcl}
n \deg V - r \deg E & = & (n \deg f^* V - r \deg f^* E)
(\deg f)^{-1} \\ [2mm]
& \geq & (n \deg G - r \deg F) (\deg f)^{-1} \mbox{~~~~by (3)} \\
[2mm]
& = & (- n \deg \widetilde{G} + r \deg F^* ) ( \deg f)^{-1}
\\
& \geq & 0 \mbox{~~~~by (4)}.
\end{array}
$$
This proves (2).

The proof in the stable case is similar.

\blem \label{l12}
Let $X^m$ and $Y^n$ be smooth projective varieties, and
$D_X$  and $D_Y$ ample divisors
 on $X$ and $Y$.  Let $\eta = a D_X + b D_Y$, $a,b >0$.
Suppose that $E$ is a vector bundle on $X \times Y$, such that for
generic $x \in X$, $y \in Y$, $E_x \simeq E|_{\{x\} \times Y}$ and $E_y
\simeq E |_{X \times \{y\}}$ are respectively
$D_Y$-semi-stable and
$D_X$-semi-stable.  Then $E$ is $\eta$-semi-stable.

Further, if either $E_x$ or $E_y$ is stable, then $E$ is
stable.
\elem
{\bf Proof.}
Let $F \subset E$ be a subsheaf.
Since $\sing F$ has codimension $\geq 2$, we can choose $x
\in X$ and $y \in Y$ such that $\sing F_x$ and $\Sing F_y$
also have codimension $\geq 2$.  Thus any torsion in $F_x$
or $F_y$ is supported in codimension $\geq 2$ and does not
contribute to $c_1(F_x)$ or $c_1(F_y)$.
Let  $\rank (E) = n$, $\rank (F) = r$.  Then we need to show
that
$$
\frac{\deg F}{r} \leq \frac{\deg E}{n} \eqno (5)
$$
assuming $E_x$ and $E_y$ are semi-stable.  Now,
$$
\begin{array}[b]{rcl}
\deg E & = & c_1(E) \cdot [aD_X + bD_Y]^{m+n-1} \\ [2mm]
& = & c_1(E) [\lambda D_X^m \cdot D_Y^{n-1} + \mu D_X^{m-1}\cdot D_Y^n]
\mbox{~~~~for some}~~
\lambda ,\mu > 0 \\ [2mm]
& = & [c_1(E_x) + c_{1,1} (E) + c_1(E_y)] \cdot [\lambda D_X^m \cdot
D_Y^{n-1} + \mu D_X^{m-1}\cdot D_Y^n] \\ [2mm]
& = & [c_1 (E_x) \cdot D_Y^{n-1}] \cdot \lambda D_X^m + [c_1(E_y)\cdot
D_X^{m-1}] \cdot \mu D_Y^n. \\
\end{array} \eqno (6)
$$
We have a similar expression for $\deg F$,
$$
\deg F = [c_1 (F_x) \cdot D_Y^{n-1}] \cdot \lambda D_X^m + [c_1(F_y)
\cdot D_X^{n-1}] \cdot \mu D_Y^n. \eqno (7)
$$
(5) follows trivially by comparing the terms in (6) and (7) and using
semi-stability of $E_x$ and $E_y$.  The rest of the lemma follows in a
similar fashion.

Before stating the next lemma, we recall very briefly some facts on
spectral curves.  For  details see [BNR] and [Li].

Let $K = K_C$ be the canonical bundle and let $W = \oplus_{i=1}^n
H^0(C,K^i)$.  Let $s=(s_1,\cdots ,s_n) \in W$, and let $C_s$ be the
associated spectral curve.  Then we have a morphism
$$
\pi : C_s \lr C
$$
of degree $n$, such that for $x \in C$, the fibre $\pi^{-1}(x)$ is
given by points $y \in K_x$ which are zeros of the
polynomial
$$
f(y) = y^n + s_1(x) \cdot y^{n-1}  + \cdots + s_n(x) .
$$

The condition that $x$ be unramified is that
the resultant $R(f,f')$ of $f$ and its derivative $f'$ be non-zero at
the point $(s_1(x),\cdots ,s_n(x))$.  Note that $R(f,f')$ is a
polynomial in the $s_i(x)$, $i=1,\cdots ,n$.

\blem \label{l13}
Given $x \in C$, there exists a smooth spectral curve $C_s$
such that the
covering map $\pi : C_s \lr C$ is unramified at $x$.
\elem
{\bf Proof.}  Note first that, if $x \in C$, there exists
$s = (s_1,\cdots ,s_n) \in W$ such that
$$
R(f,f') (s_1(x),\cdots ,s_n(x)) \neq 0.
$$
Indeed, since $|K^i|$ has no base points, given any
$(\alpha_1,\cdots ,\alpha_n) \in
\bigoplus_{i=1}^n K_x^i$, there exist
$s_i \in H^0(C,K^i)$ such that $s_i(x) = \alpha_i$, $
i=1,\cdots ,n$.

Observe that this is clearly an open condition on $W$.

Further, the subset $\{ s \in W ~|~ C_s$ is smooth$\}$ is a non-empty
open subset of $W$ [BNR, Remark 3.5] and the lemma follows.

Let $M_{\xi}$ and $\cal U_{\xi}$ be as in the
introduction, and let $\Theta_{\xi}$ denote the
restriction of the generalized theta divisor to $M_{\xi}$.
\bprop \label{l14}
Let $\cal U_{\xi}$ be the Poincar\'e bundle on $C \times M_{\xi}$ and $x \in
C$.  Then the bundle
$$
{\cal U}_{\xi ,x}\cong \cu_{\xi} |_{\{x\} \times M_{\xi}}
$$
is $\Theta_{\xi}$-semi-stable on $M_{\xi}$.
\eprop
{\bf Proof.}  For the point $x \in C$ above, choose a spectral curve
$C_s$ by Lemma 1.3, so that
$$
\pi : C_s \lr C
$$
is unramified at $x$.  Let $\pi^{-1} (x) = \{ y_1,\cdots ,y_n \}$,
$y_i$ being distinct points in $C_s$.

Let $J^{\delta}(C_s)$ denote the variety of line bundles
of degree $$\delta=d-\deg\pi_*(\co_{C_s})$$
on $C_s$, and let $P_s$ denote the subvariety of
$J^{\delta}(C_s)$ consisting of those line bundles $L$ for
which the vector bundle $\pi_*L$ has determinant $\xi$.
($P_s$ is a translate of the Prym variety of $\pi$.) Let
$\cl$ denote the restriction of the Poincar\'e bundle on
$C_s \times J^{\delta}(C_s)$ to $C_s\times P_s$. Then
[BNR, proof of Proposition 5.7], we have a dominant
rational map 
defined on an open subset  $T_s$ of $P_s$ such that $\codim (P_s - T_s) \geq
2$, and
$$
\phi : T_s \lr M_{\xi}
$$
is generically finite. The morphism $\phi$ on $T_s$ is defined by the family
$(\pi\times 1)_*{\cl}$ on $C_s\times T_s$; so, by the universal property of
$\cu_{\xi}$, we have
$$
(\pi \times 1)_* {\cal L} \simeq (1 \times \phi )^* {\cal U}_{\xi}\otimes
p_T^*L_0\eqno (8)
$$
for some line bundle $L_0$ on $T_s$. (Here $p_T : C_s\times T_s\lr T_s$ is the
projection.)  By $(8)$ we have
$$
\phi^* {\cal U}_{\xi ,x} \simeq [(\pi \times 1)_* {\cal L}
]_x\otimes L_0^{-1} ~{\rm on} ~T_s.
$$
But $[(\pi \times 1)_* {\cal L}]_x \simeq \oplus_{i=1}^n
{\cal
L}_{y_i}$ on $P_s$.  Hence
$$
\phi^* {\cal U}_{\xi ,x} \simeq \oplus_{i=1}^n ({\cal L}_
{y_i}\otimes L_0^{-1})  ~{\rm on}~
T_s. \eqno (9)
$$

We observe that the ${\cal L}_{y_i}\otimes L_0^{-1}$ are
the restrictions to $T_s$ of algebraically equivalent line
bundles on $P_s$.  Further one knows [Li, Theorem 4.3]
that $\phi^* \Theta_{\xi}$ is a multiple of the
restriction of the usual theta divisor on $J(C_s)$ to
$T_s$.

Now we are in the setting of Lemma 1.1 and we conclude that ${\cal
U}_{\xi ,x}$ is semi-stable with respect to $\Theta_{\xi}$ on $M_{\xi}$.

\bth \label{th15}
The Poincar\'e bundle $\cu_{\xi}$ on $C \times M_{\xi}$ is stable with respect
to any polarisation.
\eeth
{\bf Proof.}  Since $\hbox{Pic}\, M_{\xi} = {\Bbb Z}$,
$\hbox{Pic}\, (C \times M_{\xi}) = \hbox{Pic}\,
C \oplus \hbox{Pic}\, M_{\xi}$.  Thus, any polarisation
$\eta$ on $C \times M_{\xi}$
can be expressed in the form
$$\eta = a \alpha + b \Theta_{\xi}, ~~~ a,b > 0.$$
for some ample divisor $\alpha$ on $C$.

By Proposition 1.4, $\cu_{\xi ,x}$ is semi-stable with respect to
$\Theta_{\xi}$ for all $x \in C$ and by definition $\cu
_{\xi}|_{C\times\{m\}}$ is stable with respect to any
polarisation on $C$.  Hence by Lemma 1.2, $\cu_{\xi}$ is
stable with respect to $\eta$ on $C \times M_{\xi}$.

Note that Proposition 1.4  remains true if we replace
$M_{\xi}$, $\cu_{\xi}$ and $\Theta_{\xi}$ by $M$, $\cu$ and
$\Theta$.  (The key point is that [Li, Theorem 4.3] is valid
in this context).  We deduce at once

\bth
$\cu$ is stable with respect to any polarisation of the form
$$a \alpha  + b \Theta ,~~ a,b > 0,$$
where $\alpha$ is ample on $C$ and $\Theta$ is the
generalized theta divisor on $M$.
\eeth

\brem \label{r16}
Since $C \times M$ is a K\"ahler manifold, then by a theorem of
Donaldson-Uhlenbeck-Yau, $\cu$ admits an Hermitian-Einstein metric.
One can expect that the restriction of this metric to each factor
is precisely the metric on the factor.  It would be interesting to
know an explicit description of the metric on $\cu$.  Note
that [Ke2] contains such a description for the Picard
sheaf  in the case $g=1$, $n=1$.
\erem

\renewcommand{\thesection}{\S \arabic{section}}
\section{Some properties of End $\cu_{\xi}$}
\renewcommand{\thesection}{\arabic{section}}

Our first object in this section is to calculate the
dimension of some of the cohomology spaces of $\eee
\cu_{\xi}$.  We denote by $p: C \times M_{\xi} \lr M_{\xi}$
and $p_C: C \times M_{\xi} \lr C$ the projections.

\bprop
Let $h^i (\eee \cu_{\xi}) = \dim H^i (\eee \cu_{\xi})$.
Then,
$$h^0(\eee \cu_{\xi} ) = 1, ~~ h^1(\eee \cu_{\xi}) =g, ~~h^2(\eee \cu_{\xi}) =
3g-3.$$
\eprop
{\bf Proof.}  We can write $\eee \cu_{\xi} \cong \co \oplus
\ad \cu_{\xi}$; hence
$$H^i(C \times M_{\xi}, \eee \cu_{\xi}) = H^i(C\times
M_{\xi}, \co ) \oplus H^i (C \times M_{\xi}, \ad \cu_{\xi}
).$$
Since $\cu_{\xi} \lm_{C \times \{ m\}}$ is always stable,
we have $R_p^0(\ad \cu_{\xi}) = 0$.  So by the Leray
spectral sequence and the fact that $R_p^1(\ad \cu_{\xi})$
is the tangent bundle $T M_{\xi}$ of $M_{\xi}$, we have
$$\begin{array}{rl}
 H^i (C \times M_{\xi},\ad \cu_{\xi} ) & \cong H^{i-1}
(M_{\xi}, R_p^1(\ad \cu_{\xi}))  \\ [2mm]
& \cong H^{i-1} (M_{\xi},TM_{\xi})  \\
\end{array}$$
By [NR, Theorem 1], this space is 0 if $i \neq 2$, and has
dimension $3g-3$ if $i=2$.

On the other hand, since $M_{\xi}$ is unirational, it follows
from the K\"unneth formula that
$$H^i (C \times M_{\xi} , \co ) = H^i (C,\co_C ).$$
This space has dimension 1 if $i=0, ~g$ if $i=1$ and 0
otherwise.  The proposition follows.

\noindent
\brem \label{r22}   In fact the proof shows that
$$h^i (\eee \cu_{\xi} ) = 0~~{\rm if}~~ i > 2.$$
\erem

\blem \label{l23}
Let $L \in J (C)$ and suppose that
$E \cong E \otimes L$ for all $E \in M_{\xi}$.
Then $\cu_{\xi} \cong \cu_{\xi} \otimes p_C^* L$.
\elem
{\bf Proof.}  Since $E \cong E \otimes L$ and $\cu_{\xi} \otimes p_C^* L$ is a
family of stable bundles, there is a line bundle $L_1$
over $M_{\xi}$ such that
$$\cu_{\xi} \cong \cu_{\xi} \otimes p_C^*L \otimes p^*
L_1.$$ Fix $x \in C$, then $\cu_{\xi ,x} \cong \cu_{\xi ,x}
\otimes L_1$ over $M_{\xi}$.  Hence $c_1 (\cu_{\xi ,x})
= c_1 (\cu_{\xi ,x}) + nc_1 (L_1 )$ so
$nc_1(L_1) = 0$.  But $\hbox{Pic}\, \cm_{\xi} \cong {\Bbb
Z}$ (see [DN]); so $c_1 (L_1)=0$, which implies that $L_1$
is the trivial bundle.

The next lemma will also be required in \S 3.
\blem \label{l24}
If $\cu_{\xi} \cong \cu_{\xi} \otimes p_C^*L$ then $L \cong
\co_C$. \elem
{\bf Proof.}  If $\cu_{\xi} \cong \cu_{\xi} \otimes p_C^*L$ then
$$\begin{array}{rcl}
\co \oplus \ad \cu_{\xi} & \cong & \eee \cu_{\xi} \\ [2mm]
& \cong & \eee \cu_{\xi} \otimes p_C^* L \\ [2mm]
& \cong & p_C^* L \oplus \ad \cu_{\xi} \otimes p_C^* L. \\
\end{array}$$
Hence $H^0(C \times M_{\xi}, p_C^*L)$ and $H^0(C \times
M_{\xi}, \ad \cu_{\xi} \otimes p_C^*L)$ cannot both be
zero.  Suppose there is a non-zero section $\phi: \co \lr \ad \cu_{\xi}
\otimes p_C^* L$.
For some $x \in C$, the restriction of $\phi$ to  $\{ x\} \times
M_{\xi}$ will define a non-zero section of $\ad
\cu_{\xi,x}$, which is a contradiction since $H^0(M_{\xi},
\ad \cu_{\xi ,x})=0$ (see [NR, Theorem 2]).  Hence $H^0(C
\times M_{\xi}, \ad \cu_{\xi} \otimes p_C^* L) = 0$.

Therefore $H^0(C \times M_{\xi}, p_C^*L) \ne0$. Since
$\deg L=0$, this implies $L \cong \co_C$.

\brem
The proof of Lemma 2.4 fails when $g=1$ since [NR,
Theorem 2] is not then valid. In fact, Lemma 2.4 and the
remaining results of this section are false for $g=1$.
\erem
 We
show next that a general stable bundle $E$ is not
isomorphic to $E\otimes L$ unless $L\cong\co_C$.

\bprop \label{p25}
There exists a proper closed subvariety $S$ of $M_{\xi}$
such that, if $E\not\in S$, then $$E\cong E\otimes
L\Longrightarrow L\cong\co_C.$$ \eprop
{\bf Proof.} For any $L$, the subset $S_L = \{ {E \in M_{\xi} |E \cong E
\otimes L}\}$ is a closed subvariety of $M_{\xi}$.
If $L\not\cong\co_C$, then, by Lemmas 2.3 and 2.4, $S_L$
is a proper subvariety. On the other hand, $S_L$ can only
be non-empty if $L^n\cong\co_C$; so only finitely many of
the $S_L$ are non-empty. Since $M_{\xi}$ is irreducible,
the union $S= \bigcup\{S_L | L\not\cong\co_C\}$is a proper
subvariety of $M_{\xi}$ as required.

\brem
It follows at once from Proposition 2.6 that the action
of $J(C)$ on $M$ defined by $E\longmapsto E\otimes L$ is
faithful. Another proof of this fact has been given in
[Li, Theorem 1.2 and Proposition 1.6]. As the following
proposition shows, our set $S$ is analogous to the set $S$
of [Li, Theorem 1.2].
\erem
\bprop \label{p27}
Let $S$ be as above and $E\in M_{\xi}$. Then $E \in S$ if and only if $\ad E$
has a line sub-bundle of
degree zero.
\eprop
{\bf Proof.}  The trivial bundle cannot be a subbundle of
$\ad E$.

If $L \in J(C)$ is a subbundle of $\ad E$, so is it  of
$\eee E$,
therefore
$$H^0(C, \eee E \otimes L^*) \neq 0.$$  Hence, there is a
non-zero map $\phi : E \otimes L \lr E$, which is an
isomorphism since $E \otimes L$ and $E$ are stable bundles
of the same slope. Hence $E\in S$.

Conversely, suppose $E \cong E \otimes L$ with $L \not
\cong \co_C$.
The isomorphism $\co_C \oplus \ad E \cong L \oplus \ad E \otimes L$
implies that $\ad E \otimes L$ has a section, i.e.\ there
is a non-zero map $\phi: L^* \lr \ad E$.  Since $L^*$ and
$\ad E$ have the same slope and $\ad E$ is semi-stable
(being a subbundle of a semi-stable bundle $\eee E$ with
the same slope), $\phi$ is an inclusion.

\bcor
If $E$ is a general stable bundle of rank $2$ and determinant $\xi$, then $\ad
E$ is stable.
\ecor
{\bf Proof.}
Note that $\ad E$ has rank 3 and degree $0$ and is
semi-stable.  By the Proposition, $\ad E$ has no line
subbundle of degree 0.  On the other hand $\ad E$ is
self-dual, so it cannot have a quotient line bundle of
degree 0.
\brem \label{r38}
It would be interesting to know if $\ad E$ is
stable for a general stable
bundle $E$ of rank greater than $2$.  It is certainly true
that $\ad E$ is semistable and also that it is stable as an
orthogonal bundle [R].
\erem

\bth
If $n=2$, $\ad \cu_{\xi}$ is stable with respect to any
polarisation on $C \times M_{\xi}$.
\eeth
{\bf Proof.}  In view of Corollary 2.9 and Lemma 1.2, we
need only prove that $\ad \cu_{\xi ,x}$ is semi-stable for
some $x \in C$.  The argument is the same as for Proposition
1.4; indeed (9) shows at once that $\phi^* \ad \cu_{\xi ,x}$
can be expressed as a direct sum of restrictions to $T_s$ of
algebraically equivalent line bundles on $P_s$.

For $n > 2$, we can show similarly that $\ad \cu_{\xi}$ is
semi-stable.

\renewcommand{\thesection}{\S \arabic{section}}
\section{Deformations}
\renewcommand{\thesection}{\arabic{section}}

As in the introduction, let $H$ be any ample divisor on $C
\times M_{\xi}$, let $M(\cu_{\xi})$ denote the moduli
space of $H$-stable bundles with the same numerical
invariants as $\cu_{\xi}$ on $C \times M_{\xi}$, and let
$M(\cu_{\xi})_0$ denote the connected component of
$M(\cu_{\xi})$ which contains $\cu_{\xi}$.  One can define a
morphism
$$\beta : J (C) \lr M(\cu_{\xi})_0$$
by
$$\beta (L) = \cu_{\xi} \otimes p_C^*L.$$
Our object in this section is to prove

\bth
$\beta$ is an isomorphism.
\eeth
\brem
Note that this implies in
particular that $M(\cu_{\xi})_0$ is independent of the
choice of $H$, and is a smooth projective variety of
dimension $g$. Since $h^2(\eee \cu_{\xi}) \neq 0$, there
is no a priori reason why this should be so. \erem
{\bf Proof of Theorem 3.1.}
By Lemma 2.4, $\beta$ is injective.  Moreover the Zariski
tangent space to $M(\cu_{\xi})_0$ at $\cu_{\xi} \otimes
p_C^*L$ can be identified with $H^1(C \times M_{\xi}, \eee
\cu_{\xi})$, which has dimension $g$ by Proposition 2.1.
It
follows that, at any point of ${\rm Im }~\beta$,
$M(\cu_{\xi})_0$ has dimension precisely $g$ and is smooth.
Hence by Zariski's Main Theorem, $\beta$ is an open
immersion.  Since $J(C)$ is complete, it follows that
$\beta$ is an isomorphism.

\vspace{.15in}\noindent
V. Balaji\\*School of Mathematics\\*SPIC Science
Foundation\\*92 G. N. Chetty Road\\*T. Nagar\\*Madras 600
017\\*India\\*balaji@ssf.ernet.in
\\[.15in]
L. Brambila Paz\\*Departamento de Matematicas\\*UAM -
Iztapalapa\\*Mexico, D. F.\\*C. P.
09340\\*Mexico\\*lebp@xanum.uam.mx
\\[.15in]
P. E. Newstead\\*Department of Pure Mathematics\\*The
University of Liverpool\\*Liverpool\\*L69
3BX\\*England\\*newstead@liverpool.ac.uk

\end{document}